\documentclass[11pt,twoside]{amsart}
\usepackage{amsmath, amsthm, amscd, amsfonts, amssymb, graphicx, color}
\usepackage[bookmarksnumbered, plainpages]{hyperref}
\newenvironment{thm}{\subsection{}{\textbf {Theorem.}}\em}{}
\newenvironment{prop}{\subsection{}{\textbf {Proposition.}}\em}{}
\newenvironment{cor}{\subsection{}{\textbf {Corollary.}}\em}{}
\newenvironment{lem}{\subsection{}{\textbf {Lemma.}}\em}{}

\newenvironment{pf}{\noindent{\em
{Proof.}}}{\begin{flushright}\eop \end{flushright}\smallskip}

\newenvironment{defn}{\subsection{}{\textbf 
{Definition.}}\em}{\smallskip}

\newenvironment{rem}{\subsection{}{\textbf {Remark.}}}{\smallskip}

\newcommand\cA{\ensuremath{\mathcal A}}

\newcommand\eop{{{\hfil \ensuremath \Box}}}
\begin{document}
%------------------------------------------------------------------------------------%
%%Don not change any thing in this part
\hskip -0.2 cm
\begin{tabular}{c r}
\vspace{-0.6cm}
\href{http://jims.iranjournals.ir}{\scriptsize  \rm http://jims.iranjournals.ir}\\
\end{tabular}
\hskip 2 cm
\begin{tabular}{l}
\hline
\vspace{-0.2cm}
\scriptsize \rm\bf Journal of the Iranian Mathematical Society\\
\vspace{-0.2cm}
\scriptsize \rm ISSN (print): ????-????, ISSN (on-line): ????-???? \\
\vspace{-0.2cm}
\scriptsize Vol. {\bf\rm x} No. x {\rm(}201x{\rm)}, pp. xx-xx.\\
\scriptsize $\copyright$ 201x Iranian Mathematical Society\\
\hline
\end{tabular}
\hskip 2 cm
\begin{tabular}{c c}
\vspace{-0.1cm}
\href{http://en.ims.ir/}{\scriptsize \rm http://en.ims.ir/}\\
\vspace{-1cm}
\end{tabular}
\vspace{1.3 cm}

%------------------------------------------------------------------------------------%

\title{Arens regularity of ideals in $A(G)$, $A_{cb}(G)$ and $A_M(G)$}
\author{Brian Forrest, John Sawatzky and  Aasaimani Thamizhazhagan
 }

\thanks{{\scriptsize
\hskip -0.4 true cm MSC(2010): Primary: 43A07; Secondary: 43A22, 46J10, 47L25
Received: dd mmmm yyyy, Accepted: dd mmmm yyyy.\\
This research is partially supported by the Natural Sciences and Engineering Research Council of Canada.}}
\maketitle

\begin{abstract}  In this paper, we look at the question of when various ideals in the Fourier algebra $A(G)$ or its closures $A_M(G)$  and $A_{cb}(G)$ in, respectively, its multiplier and $cb$-multiplier algebra are Arens regular. We show that in each case, if a non-zero ideal is Arens regular, then the underlying group $G$ must be discrete. In addition, we show that if an ideal $I$ in $A(G)$ has a bounded approximate identity, then it is Arens regular if and only if it is finite-dimensional. 
\end{abstract}

\vskip 0.2 true cm

%------------------------------------------------------------------------------------%

\pagestyle{myheadings}
\markboth{\rightline {\sl  J. Iranian Math. Soc. x no. x (201x) xx-xx \hskip 4 cm  B. Forrest, J. Sawatzky and A. Thamizhazhagan }}
         {\leftline{\sl J. Iranian Math. Soc. x no. x (201x) xx-xx \hskip 4 cm  B. Forrest, J. Sawatzky and A. Thamizhazhagan}}

\bigskip
\bigskip

%------------------------------------------------------------------------------------%
%------------------------------------------------------------------------------------%

\section{Introduction}

Let $G$ be a locally compact group with a fixed left Haar measure $\mu$. Let $\Sigma_G$ denote the collection of equivalence classes of weakly continuous unitary representations of $G$.  The Fourier-Stieltjes algebra, $B(G)$ is the space of all coefficient functions of weakly continuous unitary representations on $G$. That is \[B(G)=\{u(x)=\langle\pi(x)\xi,\eta\rangle| \pi\in \Sigma_G, \xi, \eta \in \mathcal{H}_{\pi}\}.\]
$B(G)$ is a commutative Banach algebra under pointwise operations when given the norm it inherits as the dual of the group $C^*$-algebra $C^*(G)$. 

The left regular representation $\lambda$ on $L^2(G)$  is defined
\[\lambda(x)(f)(y)=f(x^{-1}y)\]
for each $x,y \in G$. The Fourier algebra $A(G)$ consists of all coefficient functions of $\lambda$. That is all functions of the  form 
\[u(x)=\langle\lambda(x)f,g\rangle\] where $f, g \in L^2(G)$. The Fourier algebra is a closed ideal of $B(G)$ with $\Delta(A(G))=G$. It is also the predual of the von Neumann subalgebra $VN(G)$ of $B(L^2(G))$ generated by 
$\{\lambda(x)\mid x\in G\}$. It can also be realized as the closure of
the space of elements in $B(G)$ with compact support. (See \cite{Eymard} 
for more on the nature of $B(G)$ and $A(G)$.)

It is well known that the multiplication on any Banach algebra $\cA$ can be extended to its second dual in two natural ways as demonstrated by Arens in \cite{Arens}. Generally, these two Arens products are different. If they agree, we say that the Banach algebra $\cA$ is \textit{Arens regular}. For a commutative Banach algebra $\cA$, Arens regularity occurs precisely when the second dual $\cA^{**}$ is commutative. 

The most well-known examples of Arens regular algebra are the closed subalgebras of $B(\mathcal{H})$, the bounded operators on a Hilbert space $\mathcal{H}$.  Typically, many of the algebras arising from locally compact groups fail to be Arens regular unless the group structure is significantly restricted. For example, the first author showed that if $A(G)$ is Arens regular, then $G$ must be discrete and every amenable subgroup of $G$ must be finite \cite[Theorem 3.2 and Proposition 3.7]{For2} extending an earlier result of Lau's for amenable groups \cite[Proposition 3.3]{Lau2}. In fact, it is conjectured, and strongly believed to be true, that whenever $A(G)$ is Arens regular, $G$ must be finite.  

The literature concerning Arens regularity of algebras arising from locally compact groups is extensive. For the Fourier algebra in particular see   for example \cite{For1}, \cite{For2},  \cite{Hu}, \cite{Lau2} and  \cite{Ulger2}

In this paper, we will consider the possible Arens regularity of ideals in $A(G)$ and in two related Banach algebras $A_M(G)$ and $A_{cb}(G)$, that arise from $A(G)$ as its closure in its space of multipliers and completely bounded multipliers respectively. We will show that if $I$ is a non-zero closed ideal in any of these three algebras that is Arens regular, then $G$ must be discrete.

\bigskip

\bigskip

\section{Preliminaries and Notation} 

Throughout this paper, $\cA$ will denote a  Banach algebra. In this case, the dual $\cA^*$ becomes a Banach $\cA$-bimodule with respect to the module actions 
\[\langle u\cdot T, v\rangle =\langle T, vu\rangle  \]
and 
\[\langle  T\Box u, v\rangle  = \langle T, uv\rangle  \]
for every $u,v\in \cA$ and $T \in \cA^*$.

It is well known that there are two natural products that can be used to extend the multiplication of $\cA$ to its second dual $\cA^{**}$. These two Arens products are  defined as follows: 
\begin{itemize}
\item[1a)] $\langle u\cdot T, v\rangle  =\langle T, vu\rangle $ for every $u,v \in \cA$ and $T \in \cA^*$.
\item[1b)] $\langle T\odot m, u\rangle  = \langle m, u \cdot T\rangle $ for every $u \in \cA$ and $T \in \cA^*$ and $m\in \cA^{**}$.
\item[1c)] $\langle m\odot n, T\rangle  = \langle n, T\odot m\rangle $ for every $T \in \cA^*$ and $m,n \in \cA^{**}$.
\item[2a)] $\langle T\Box u, v\rangle  =\langle T, uv\rangle $ for every $u,v \in \cA$ and $T \in \cA^*$.
\item[2b)] $\langle m\Box T, u\rangle  = \langle m, T\Box u\rangle $ for every $u \in \cA$ and $T \in \cA^*$ and $m\in \cA^{**}$.
\item[2c)] $\langle m\Box n, T\rangle  = \langle m, n\Box T\rangle $ for every $T \in \cA^*$ and $m,n \in \cA^{**}$.
\end{itemize} 

\begin{rem} If $\cA$ is a commutative Banach algebra, then it is easy to see that $u\cdot T=T\Box u$ for every $u \in \cA$ and $T \in \cA^*$. In particular, if $m\in \cA^{**}$ and $u\in \cA$, then 
\[m\odot u=u\odot m= m\Box u = u\Box m.\]
However, as is well known, even if $\cA$ is commutative this does not mean that the two multiplications agree on $\cA^{**}$. Moreover, even though $\cA$ is assumed to be commutative, it is generally not the case that $\cA^{**}$ would be as well. 
\end{rem} 
\begin{defn}
If $\cA$ is a Banach algebra for which $m\odot n=m\Box n$ for every $m,n \in \cA^{**}$, we say that $\cA$ is \textit{Arens regular}. 
\end{defn}

 From here on we will assume that $\cA$ is a commutative Banach Algebra with maximal ideal space $\Delta(\cA)$.  Moreover, if we are speaking of $\cA^{**}$, we will assume that the product we are using is $\odot$ unless otherwise specified. 

 We will proceed with the following definitions and notational conventions. 

\begin{defn} We call the space 
\[UCB(\cA)=\overline{\textnormal{span}\{v \cdot T\mid v\in \mathcal{\cA}, T\in \cA^* \}}^{-\|\cdot\|_{\cA^*}}\]  
the \textit{uniformly continuous functionals on} $\cA$. 

We call $T\in \cA^*$ a \textit{(weakly) almost periodic functional} on $\cA$ if 
\[\{ u \cdot T\mid u \in \cA, \|u\|_{\cA} \leq 1\}\]
is relatively (weakly) compact in $\cA^*$ and we denote the space of all (weakly) almost periodic functionals on  $\cA$ by $AP(\cA)$ ($WAP(\cA)$).

\end{defn} 

\begin{rem} It is a well-known criterion of Grothendieck that $T\in \cA^*$  is weakly almost periodic if and only if given two nets $\{u_{\alpha}\}_{\alpha \in \Omega_1}$ and  $\{v_{\beta}\}_{\beta\in \Omega_2}$ in $\cA$ we have that 
\[\lim\limits_{\alpha} \lim\limits_{\beta}\langle T,u_{\alpha}v_{\beta}\rangle=
\lim\limits_{\beta} \lim\limits_{\alpha}\langle T,u_{\alpha}v_{\beta}\rangle\]
whenever both limits exist. From here it is easy to see that $\cA$ is Arens regular if and only if $WAP(\cA)=\cA^*$. Moreover, as a consequence, this shows that if $\cA$ is commutative, then $\cA$ is Arens regular if and only if $\cA^{**}$ is also commutative in either, and hence both Arens products.(See \cite{DuncHoss} and \cite{Pym}.)

\end{rem}

\begin{defn} We say that a closed subspace $X\subseteq \cA^*$ is \textit{invariant} if $u \cdot T \in X$ for every $u\in \cA$ and $T\in X$. \vspace{.05in}

Given a closed invariant  subspace $X$ of $\cA^*$ and an $m\in X^*$ , we define the linear operator $m_L:X\to \cA^*$ by 
\[\langle m_L(T),u\rangle :=\langle m,u\cdot T \rangle \]
for every $T\in \cA$ and $u\in \cA$. We say that $X$ is topologically introverted if $m_L(T)\in X$ for every $m\in X^*$ and $T \in X$. 

It is also well known that if $X$ is topologically introverted then $X^*$ can be made into a Banach algebra by mimicking what we did for $\cA^{**}$ as follows. 
\begin{itemize}
\item[1)] For each $T\in X$ and $m\in X^*$, we define $T\odot m = m_L(T)$. 
\item[2)] For each $T\in X$ and $n,m\in X^*$, we define $\langle m\odot n,T\rangle  = \langle n,  T\odot m\rangle $
\end{itemize}

\end{defn} 
It is well-known and straightforward to show that all three of  $AP(\cA),  WAP(\cA)$, and $UCB(\cA)$ are closed introverted subspaces of $\cA^*$. Moreover, if $\cA$ is commutative and $X$ is topologically introverted, then $X^*$ is commutative if and only if $X\subseteq WAP(\cA)$. 

%(See for example \cite[Proposition 3.1 and Proposition 3.2]{Lau}).

\section{Multipliers of the Fourier Algebra} 

Let $G$ be a locally compact group. We let $A(G)$ and $B(G)$ denote the Fourier and
 Fourier-Stieltjes algebras of $G$, which are Banach algebras of continuous functions
on $G$ and were introduced in [5]. 
 A \emph{multiplier} of $A(G)$ 
is a bounded continuous  function 
$v \!: G \to \mathbb{C}$
such that $v A(G) \subseteq A(G)$. Each multiplier $v$ of $A(G)$ determines a 
 linear operator $M_v$ on $A(G)$ defined by 
$M_{v}(u)=vu$ for each $u\in A(G)$. It is a routine consequence of the Closed 
Graph Theorem  that each $M_v$ is also bounded. The {\it multiplier algebra\/} 
of $A(G)$ is the closed subalgebra

\[
  MA(G) := \{ M_v : \textnormal{$v$ is a multiplier of $A(G)$} \}
\]
of $B (A(G))$, where  $B (A(G))$ denotes the algebra of all bounded linear 
operators from $A(G)$ to $A(G)$. 
Throughout this paper, we will typically use $v$ in place of the 
operator $M_{v}$ for notational convenience and we will
write $\parallel v\parallel _{MA(G)}$ to represent the norm of $M_{v}$ in  $B (A(G))$. 

Let $G$ be a locally compact group and let $VN(G)$ denote its group 
von Neumann algebra. The duality
\[
  A(G) = VN(G)_*\]
equips $A(G)$ with a natural operator space structure. Given this 
operator space structures,  we can define the 
$cb$\emph{-multiplier algebra} of $A(G)$ to be
\[
  {M}_{cb}A(G) := CB(A(G)) \cap {M}(A(G)),
\]
where $CB(A(G))$ denotes the algebra of all completely bounded 
linear maps from $A(G)$ into itself. 
We let $\parallel  v\parallel _{cb}$ denote the $cb$-norm of 
the operator $M_{v}$. 
It is well known that ${M}_{cb}(A(G))$ is a closed subalgebra 
of $CB(A(G))$ and 
is thus a (completely contractive) Banach algebra
with respect to the norm $\parallel \cdot \parallel _{cb}$.

It is clear that, 
\[A(G)\subseteq B(G) \subseteq {M}_{cb}(A(G)) \subseteq {M}(A(G)).\] 
Moreover, for  $v \in A(G)$ we have that 
\[\parallel v\parallel _{A(G)}=\parallel v\parallel _{B(G)}\ge 
\parallel v\parallel _{cb}\ge \parallel v\parallel _{M}.\]

In case $G$ is an amenable group, we have 
\[B(G)={M}_{cb}(A(G)) = {M}(A(G)) \] 
and that 
\[\parallel v\parallel _{B(G)} = \parallel v\parallel _{cb} = \parallel v\parallel _{M}\]
for any $v\in B(G)$.

\begin{defn}  Given a locally compact group $G$ let 
\[A_{M}(G)\stackrel{def}{=}A(G)^{-\parallel \cdot \parallel _{M } }\subseteq M(A(G)).\] 
and
\[A_{cb}(G)\stackrel{def}{=}A(G)^{-\parallel \cdot \parallel _{cb } }\subseteq M_{cb}(A(G)).\] 

\end{defn} 
 We refer the reader to \cite{Haag} for the basic properties of ${M}_{cb}(A(G))$.

\begin{rem} The algebra $A_{cb}(G)$ was introduced by the first author in \cite{For} where it was denoted by $A_0(G)$. In that paper, we show that in the case of 
$\mathbb{F}_2$, the free group on two generators, $A_{cb}(G)$ shares many of the properties characteristic of the Fourier-algebra of an amenable group. In particular, the algebra 
$A_{cb}(\mathbb{F}_2)$ is known to have a bounded approximate identity. The  locally compact groups $G$ for which $A_{cb}(G)$ has a bounded approximate identity are called \textit{weakly amenable groups}. All amenable groups are weakly amenable, but many classical non-amenable groups such as $\mathbb{F}_2$ and $SL(2,\mathbb{R})$ are weakly amenable. 

\end{rem}

\begin{rem} Let $\cA(G)$ denote either $A_{cb}(G)$ or $A_M(G)$. Consider the following map and its adjoints:

\smallskip
\quad  $i_{\cA}:  A(G)\rightarrow \cA(G)$    
 
\smallskip
\quad  $i_{\cA}^*:  \cA(G)^*\rightarrow VN(G)$     

\smallskip
\quad  $i_{\cA}^{**}:  VN(G)^*\rightarrow \cA(G)^{**}$,  

\noindent
where $i_{\cA}$ denotes the inclusion map. Since $i_{\cA}$ has dense range, $i_{\cA}^*$ is injective and as such is invertible
with inverse ${i_{\cA}^*}^{-1}$ on $Range (i_{\cA}^*)$. It is easy to see that $i_{\cA}^{*}$ is simply the restriction map. That is 
\[i_{\cA}^{*}(T)=T_{|_{A(G)}}.\]

It will also be useful to view all of the above maps as embeddings. That is, when $G$ is  
non-amenable $\cA(G)^*$  can be viewed as a proper subset of $VN(G)$ and $VN(G)^{*}$ as a proper subset of $\cA(G)^{**}$. 

\end{rem}

\section{Topological Invariant Means} 

\begin{defn}  Let $\cA$ be a commutative Banach algebra with maximal ideal space $\Delta(\cA)$. Let $X$ be a closed submodule of $\cA^*$ containing $\phi \in \Delta(\cA)$.  Then $m\in X^*$ is called a topologically invariant mean (TIM) on $X$ at $\phi$ if
\begin{itemize} 
\item[i)] $\| m\|_{X^{*}} = \langle m, \phi \rangle  = 1 $,
\item[ii)] $ \langle m,  v\cdot T \rangle  =  \phi(v) \langle m,    T \rangle  $ for every $v\in \cA$ and $T\in X$.
\end{itemize}

We denote the set of topologically invariant means on $X$ at $\phi$ by $TIM_{\cA}(X,\phi)$. 

\end{defn} \vspace{.05in}

\textbf{Note:} For the rest of this section we will focus our attention on the algebras $A(G)$, $A_{cb}(G)$, $A_M(G)$, and their closed ideals.

\begin{defn} Let $\cA$ be one of the algebras  $A(G)$, $A_{cb}(G)$ or $A_M(G)$. Let $x\in G$ define the isometry $L_x:\cA \to \cA$ by 
\[L_x(u)(y) = u(xy),\]
for each $y \in G$. 

\end{defn}

The following proposition will prove useful. 

\begin{prop}\label{UCB}  Let $\cA(G)$ be one of the algebras  $A(G)$, $A_{cb}(G)$ or $A_M(G)$. Let $X\subseteq \cA(G)^*$ be a closed submodule. Let $x\in G$. 
\begin{itemize} 
\item[i)] If $Y= L^*_x(X)$, then $Y$ is a closed submodule of $\cA(G)^*$ and
\[u\cdot L^*_x(T)=L_x^*(L_x(u) \cdot T)\]
for every $u\in \cA$ and $T \in X$. 

\item[ii)] Let $x\in G$. Then $L^*_x(\phi_e)=\phi_x$ where $e$ denotes the identity of $G$. 
\item[iii)] Let $m\in TIM_{\cA(G)}(X,\phi_e)$. If   $x\in G$, Then $\phi_x\in L_x^*(X)$ and
 $L^{**}_{x^{-1}}(m) \in TIM_{\cA(G)}(L^*_x(X),\phi_x)$. 
\end{itemize}
\end{prop} 

\begin{pf} Let $\cA(G)$ be one of the algebras  $A(G)$, $A_{cb}(G)$ or $A_M(G)$. Let $X\subset \cA(G)^*$ be a closed submodule. Let $x\in G$. 
\begin{itemize} 
\item[i)] Since $L_x^*$ is an isometry, it is clear that $Y$ is a closed subspace. 

Let $T\in X$ and $u,v\in \cA(G)$ . Observe that  

\begin{eqnarray*}
\langle u\cdot L^*_x(T), v \rangle  &=& \langle   L^*_x(T) ,uv\rangle \\
&=& \langle  T,  L_x( uv)\rangle \\
&=&\langle T,L_x(u)L_x(v)_\rangle \\
&=& \langle L_x(u) \cdot T, L_x(v)\rangle \\
&=& \langle L^*_x(L_x(u) \cdot T),v\rangle \\
\end{eqnarray*}
 Hence $u\cdot L^*_x(T)=L_x^*(L_x(u) \cdot T) \in L^*_x(X)=Y$. 

\item[ii)] Let $u\in \cA(G)$. Then 
\[\langle L^*_x(\phi_e),u\rangle =\langle \phi_e, L_x(u)\rangle =L_x(u)(e)=u(x)=\langle \phi_x,u\rangle .\]
This shows that $L^*_x(\phi_e)=\phi_x$. 

\item[iii)]  Let $T=L^*_x(T_1)$ with $T_1\in X$. If $u\in \cA(G)$, we have 
\begin{eqnarray*}
\langle L^{**}_{x^{-1}}(m), u\cdot T \rangle &=&\langle m,L^*_{x^{-1}}(u\cdot T) \rangle \\
&=&\langle m,L^*_{x^{-1}}(u\cdot L^*_x(T_1)) \rangle \\
&=&\langle m,,L^*_{x^{-1}}(L_x^*((L_x(u) \cdot T_1)\rangle \\
&=&\langle m,L_x(u) \cdot T_1\rangle \\
&=&L_x(u)(e) \langle m, T_1\rangle \\
&=&u(e) \langle m, T_1\rangle \\
&=&u(e) \langle m,L^*_{x^{-1}}(L^*_{x} (T_1))\rangle \\
&=&u(e) \langle L^{**}_{x^{-1}}(m),T\rangle 
\end{eqnarray*} 

Since $L^{**}_{x^{-1}}$ is an isometry we have that 
\[\|L^{**}_{x^{-1}}(m)\|_{\cA^*}=\|m\|_{\cA^*}=1.\] Finally, we have that 
\[\langle L^{**}_{x^{-1}}(m), \phi_x\rangle =\langle m,L^{*}_{x^{-1}}(\phi_x)\rangle =\langle m,\phi_e\rangle =1.\]

\end{itemize}

\end{pf} 

The next result follows immediately from Proposition \ref{UCB} iii).

\begin{cor} Let $\cA(G)$ be any of $A(G)$, $A_{cb}(G)$ or $A_M(G)$. Let $x\in G$. Then
\[ | TIM _{\cA(G)}( \cA(G)^*, \phi_x )| = |TIM _{\cA(G)}( \cA(G)^*, \phi_e )|,\]
where $|\cdot|$ represents the cardinality of the underlying set. 

 \end{cor}

The first author together with T. Miao established  the next result  for $TIM_{A_M}(A_M^*,\phi_e)$ in \cite{ForMiao}.

\begin{prop} Let $\cA$ be either $A_{cb}(G)$ or $A_M(G)$. Let $\cA(G)\cdot VN(G) = \{ u\cdot T: u\in \cA(G), T\in VN(G) \}$. Then

\smallskip
\begin{enumerate}

\item[i)] $\cA(G)\cdot VN(G)\subseteq UCB(A(G))$

\item[ii]  $i^*(v\cdot T)=v\cdot i^*(T)$ for each $v\in \cA(G), T\in \cA(G)^{*}$.

\item[iii)] $i^{*}(UCB(\cA(G)))\subseteq UCB(A(G))$.

\item[iv)] If $\cA$ has a bounded approximate identity, then $\cA(G)\cdot VN(G)\cdot VN(G)=UCB(A(G))$ 

\item[v)] $u \cdot T \in i^{*}(\cA(G)^{*})$ for each $u\in A(G), T\in VN(G)$.

\end{enumerate}

 \end{prop}\label{UCB1}

\begin{pf} \begin{enumerate}
\item[i)] Since $UCB(A(G))$ is a closed subspace of $VN(G)$ and since $A(G)$ is dense in $\cA(G)$ with respect to the norm $\|\cdot\|_{\cA(G)}$, to establish $i)$ we need only show that for any sequence 
$\{v_{n}\}\subset  A(G)$ and $v\in \cA(G)$ with 
$\parallel v_{n} - v\parallel _{\cA(G)}\rightarrow 0$ and any $T\in VN(G)$, we have 
$\Vert v_{n} \cdot T-v \cdot T\Vert_{VN(G)}\rightarrow 0$. However, 
this follows immediately since for any $u\in A(G)$
\begin{eqnarray*}
\mid \langle v_n\cdot T-v\cdot T, u \rangle\mid &=& \mid \langle (v_n-v)\cdot T,u \rangle\mid\\
&=& \mid \langle T,(v_n-v)u \rangle\mid\\
&\leq &\parallel T\parallel _{VN(G)} 
 \parallel v_{n}-v\parallel _{\cA(G)}  \parallel u\parallel _{A(G)}.\\
\end{eqnarray*}
\item[ii)]  Let $v\in \cA(G)$,  $T\in \cA(G)^*$ and $u\in A(G)$. We have 
\begin{eqnarray*}
\langle i^{*}(v\cdot T), u\rangle &=&\langle v\cdot T, i(u)\rangle \\
&=&\langle  T, v i(u)\rangle \\
&=&\langle  T,  i(vu)\rangle \\
&=&\langle i^{*}( T), vu\rangle \\
&=&\langle v\cdot i^{*}( T), u\rangle \\
\end{eqnarray*}
Hence $i^{*}(v\cdot T)=v\cdot i^{*}( T)$. 

\item[iii)] Because of i) above we  need only show that $i^{*}( v\cdot T ) \in \cA(G)\cdot VN(G)$ for any $v\in \cA(G)$ and $T\in \cA(G)^*$. However, this follows immediately from ii).

\item[iv)] Since $\cA(G)$ has a bounded approximate identity it 
follows from Cohen's Factorization Theorem and from $i)$ above that 
$A_M(G)\cdot VN(G)$ is a closed subspace of $UCB(\hat {G})$. 
  However, since 
$A(G)\cdot VN(G)\subseteq \cA(G)\cdot VN(G)$, it is also clear that $\cA(G)\cdot VN(G)$ 
is dense in $UCB(\hat {G})$. It follows that  $\cA(G)\cdot VN(G)\cdot VN(G)=UCB(A(G))$ 
\item[v)] 
Let $u\in A(G), T\in VN(G)$. Then we can define a linear functional on $\cA(G)$ by
\[\varphi _{u,T}(v)=\langle T,uv\rangle\]
for each $v\in \cA(G)$. It is also clear that $\varphi _{u,T}$ has norm at most 
$\Vert u\Vert_{A(G)}\Vert T\Vert_{VN(G)}$. Moreover, this linear functional 
agrees with $u\cdot T$ on $A(G)$ and as such $u\cdot T=i^{*}(\varphi _{u,T})$. 

\end{enumerate}
\end{pf}

\begin{thm} Let $\cA(G)$ denote either $A_{cb}(G)$ or $A_M(G)$. For any locally compact group,  $i^{**}(TIM_{A(G)}(A(G)^*, \phi_x) \subseteq TIM_{\cA(G)} (\cA(G)^*, \phi_x) $. 
Moreover, $i^{**}: TIM_{A(G)}(A(G)^*, \phi_x)  \rightarrow TIM_{\cA(G)} (\cA(G)^*, \phi_x) $ is a bijection.
 \end{thm}

\medskip

\begin{pf}
We will first show  that $i^{**}(TIM _{A(G)}( A(G)^*, \phi_x ) ) \subseteq TIM _{\cA(G)}(\cA(G)^*, \phi_x) $.  

Let $m \in TIM _{A(G)}( A(G)^*, \phi_x )) $. Let $v \in \cA(G)$ and $T\in \cA(G)^*$.
Then there exists $\{u_n\}\subset A(G)$ such that $\parallel u_{n}-v\parallel _{\cA(G)}\rightarrow 0$.
Since $\parallel u_{n}-v\parallel _{\infty}\leq \parallel u_{n}-v\parallel _{\cA(G)}$ it follows 
that $u_{n}(x)\rightarrow v(x)$. 

Next, we note that in a similar manner to  the proof of Proposition \ref{UCB1}i), we can show that 
$  u_{n}\cdot T\rightarrow v\cdot T$ in the norm $\Vert \cdot \Vert_{\cA(G)^*}$ for each
 $T \in \cA(G)^{*}$. 
Appealing this time to Proposition \ref{UCB1}ii), it follows that
\begin{eqnarray*}
\langle i^{**}(m), v\cdot T\rangle&=&\lim\limits_{n\rightarrow \infty}
\langle i^{**}(m), u_n\cdot T\rangle\\
&=&\lim\limits_{n\rightarrow \infty}\langle m , i^*(u_n\cdot T)\rangle\\
&=&\lim\limits_{n\rightarrow \infty}\langle m ,u_n \cdot i^*(T)\rangle\\
&=&\lim\limits_{n\rightarrow \infty}u_n(x) \langle m ,  i^*(T)\rangle\\
&=&\lim\limits_{n\rightarrow \infty}v(x)\langle m , i^*(T)\rangle\\
&=&v(x)\langle i^{**}(m) , T\rangle.\\
\end{eqnarray*}
This shows that $i^{**}(TIM _{A(G)}( A(G)^*, \phi_x ))  \subseteq TIM _{\cA(G)}( \cA(G)^*, \phi_x ) $.  

We next show that $i^{**}: TIM _{A(G)}( A(G)^*, \phi_x )  \rightarrow TIM _{\cA(G)}( \cA(G)^*, \phi_x )$ is injective.
To see this, we first note that if $m_{1},m_{2}\in TIM _{A(G)}( A(G)^*, \phi_x ) $ with $m_{1}\not=m_{2}$, 
then there exists an $T\in VN(G)$ for which 
\[\langle m_{1},T\rangle\not=\langle m_{2},T\rangle.\]

Next choose $u_{0}\in A(G)$ with  $u_{0}(x)=1$. Then 
\[\langle m_{1},u_{0}\cdot T\rangle=\langle m_{1},T\rangle\not=\langle m_{2},T\rangle
=\langle m_{2},u_{0}\cdot T\rangle.\]

Since $u_{0}\cdot T\in \cA(G)^{*}$, we have 
\begin{eqnarray*}
\langle i^{**}(m_{1}),u_{0}\cdot T\rangle&=&\langle m_{1},i^*(u_{0}\cdot T)\rangle\\
&=&\langle m_{1},u_{0}\cdot T\rangle\\   
&\not=&\langle m_{2},u_{0}\cdot T\rangle\\
&=&\langle m_{2},i^*(u_{0}\cdot T)\rangle\\
&=&\langle i^{**}(m_{2}),u_{0}\cdot T\rangle\\
\end{eqnarray*}
so that $i^{**}(m_{1})\not=i^{**}(m_{2})$.

Finally, we show that $i^{**}: TIM _{A(G)}( A(G)^*, \phi_x )   \rightarrow TIM _{\cA(G)}(\cA(G)^*, \phi_x )  $ is surjective.

Let $M \in TIM _{\cA(G)}(\cA(G)^*, \phi_x ) $. 
First note that if $u,v\in A(G)$, with $u(x)=1=v(x)$ and if $T \in VN(G)$, then 
$u\cdot T$ and $v\cdot T$ are in $\cA(G)^*$ and
\[\langle M, u \cdot T\rangle =  \langle M,  v \cdot (u\cdot  T)\rangle = 
\langle M, u\cdot  (v\cdot  T) \rangle =  \langle M, v\cdot  T\rangle. \]

Pick a $u_{0}\in A(G)$ with $\parallel u_0\parallel _{A(G)}=1$ and $u_{0}(x)=1$. 
We can define $m_M\in A(G)^{**}$ by \[\langle m_M,  T\rangle  =  \langle M, u_{0} \cdot  T\rangle\]
for $T\in VN(G)$. 

Note that $\|m_M\|_{A(G)^{**}}\leq 1$. It is clear from the observation above that if $v\in A(G)$ is such that
$v(x)=1$, then $\langle m_M,  v\cdot T\rangle =\langle m_M,  T\rangle $. 
We also have that 
 \[ \langle m_M,  \phi_x\rangle  =  \langle M, u_{0} \cdot   \phi_x\rangle 
 =  \langle M, u_{0}(x) \phi_x\rangle 
 =  \langle M,  \phi_x\rangle 
 = 1. \]
 That is, $m_M\in TIM _{A(G)}( A(G)^*, \phi_x ) $.

Finally, if $T\in \cA(G)^{*}$, then 
 \[ \langle i^{**} (m_M), T\rangle = \langle m_M,  i^*(T)  \rangle  
= \langle M, u_{0}\cdot i^*(T)  \rangle 
= \langle M, u_{0}.T  \rangle 
= \langle M, T  \rangle. \]
 Therefore, $i^{**}(m_{M})=M$.

\end{pf}

\begin{defn}
Given a locally compact group $G$ we let $b(G)$ denote the smallest cardinality
of a neighbourhood basis at the identity $e$ for $G$. 

\end{defn}

The next corollary follows immediately from the previous theorem and from Hu \cite{Hu}.

\begin{cor}\label{unique} 
Let $G$ be a non-discrete
locally compact group. Let $\cA(G)$ be $A(G)$, $A_{cb}(G)$ or $A_M(G)$. Then \[\mid TIM _{\cA(G)}(\cA(G)^*, \phi_x ) ) \mid = 2^{2^{b(G)}}.\]
In particular, $\cA(G)^*$ admits a unique topological invariant mean if and 
only if $G$ is discrete. 
\end{cor}

We now turn our attention to ideals in the algebra $\cA(G)$ where $\cA(G)$ is any of  $A(G)$, $A_{cb}(G)$ or $A_M(G)$.

\begin{lem} Let $\cA(G)$ be $A(G)$, $A_{cb}(G)$ or $A_M(G)$. Let $I$ be a closed ideal in $\cA(G)$. Assume that $x\not \in Z(I)$. Let 
$M\in  TIM _{\cA(G)}(\cA(G)^*, \phi_x )$. Then $M \in (I^{\perp})^{\perp}$. 

\end{lem} 

\begin{pf} Let $T\in I^{\perp}$. Since $Z(I) $ is closed and since $x\not \in Z(I)$, we can find an open  neighborhood $U$ of $x$ such $U\cap Z(I)=\emptyset $. We can then find a $u_0 \in A(G)\cap C_{00}(G)$ so that $supp(u_0)\subseteq U$, and $u_0(x)=1$. It follows that $u_0\in I$.  

Since $T\in I^{\perp}$, we have that for any $u\in \cA(G)$ that 
\[\langle u_0 \cdot T, u\rangle  = \langle T, u_0u\rangle =0\]
so $u_0 \cdot T=0$. However, since $M\in  TIM _{\cA(G)}(\cA(G)^*, \phi_x )$, and since $u_0(x)=1$, 
\[\langle M, T\rangle =\langle M,u_0 \cdot T\rangle =0.\]
\end{pf} 

\begin{rem} The previous lemma shows that if $I$ is a closed ideal in  $\cA(G)$ with $x\not \in Z(I)$, and if $M\in  TIM _{\cA(G)}(\cA(G)^*, \phi_x )$, then we can view $M$ as an element $\hat{M}$ of $I^{**}$ in a canonical way. Specifically, if $T\in I^{**}$ and $T_1 $ is any extension of $T$ we can define 
\[\hat{M}(T)=M(T_1)\]

and $\hat{M}$ is well defined since $M\in (I^{\perp})^{\perp}$. We claim that $\hat{M}\in    TIM _{I}(I^*, \phi_{{x}_{|_{I}} })$. To see that this is the case we note that  
\[\|\hat{M}\|_{I^*}=\|M\|_{\cA(G)^{**}} = M(\phi_x )=\hat{M}(\phi_{{x}_{|_{I}} }).\] 
If $u\in I$, then if $T_{1_{|_I}}=T$, then $u\cdot T_{1_{|_I}}=u\cdot T$ and as such 
\[ \hat{M}(u\cdot T) = M(u\cdot T_1)=u(x) M(T_1)= \phi_{{x}_{|_{I}} }(u) \hat{M}(T).\]

This gives us a map $\Gamma: TIM _{\cA(G)}(\cA(G)^*, \phi_x )\to    TIM _{I}(I^*, \phi_{{x}_{|_{I}} })$, given by 
\[\Gamma{M}=\hat{M}. \]

\end{rem} The map $\Gamma: TIM _{\cA(G)}(\cA(G)^*, \phi_x )\to   TIM _{I}(I^*, \phi_{{x}_{|_{I}} })$, given by 
\[\Gamma(M)=\hat{M}. \]

\begin{thm} Let $\cA(G)$ be $A(G)$, $A_{cb}(G)$ or $A_M(G)$. Let $I$ be a closed ideal in $\cA(G)$. Assume that $x\not \in Z(I)$. The map $\Gamma: TIM _{\cA(G)}(\cA(G)^*, \phi_x )\to   TIM _{I}(I^*, \phi_{{x}_{|_{I}} })$, given by 
\[\Gamma(M)=\hat{M} \]
is a bijection.

\end{thm} 

\begin{pf} Assume that $M_1,M_2\in TIM _{\cA(G)}(\cA(G)^*, \phi_x )$ and that $T_0\in \cA(G)^*$ is such that $M_1(T_0)\not =M_2(T_0)$. Let $T\in I^*=T_{0_{|_I}}$. Then 
\[\hat{M_1}(T)=M_1(T_0)\not =M_2(T_0)=\hat{M_2}(T)\]
so 
\[\Gamma(M_1)=\hat{M_1}\not =\hat{M_2}=\Gamma(M_2)\]
and hence $\Gamma$ is injective. 

Next we let $m\in   TIM _{I}(I^*, \phi_{{x}_{|_{I}} })$. We let $u_0\in I$ be such that $u_0(x)=1$ with $\|u_0\|_{\cA(G)}=1$. First, observe that 
$(u_0\cdot T)_{|_I}=u_0\cdot (T_{|_I})$ for each  $T \in \cA(G)^*$.  Then we define $M\in \cA(G)^{**}$ by  
\[\langle M,T\rangle =\langle m,u_0\cdot (T_{|_I})\rangle =\langle m,T_{|_I}\rangle \]
for each $T \in \cA(G)^*$. We have that 
\[\|M\|_{\cA(G)^{**}}\leq \|m\|_{I^{**}} \|u_0\|_{\cA(G)}=1.\]
Moreover
\[\langle M,\phi_x\rangle = \langle m,u_0 \cdot \phi_{{x}_{|_{I}} } \rangle = u_0(x) \langle m,\phi_{{x}_{|_{I}} } \rangle =1.\]

Next let $T\in \cA(G)^*$ and let $u\in \cA(G)$. Then 
\begin{eqnarray*}
\langle M,u\cdot T\rangle &=& \langle m,u_0\cdot (u\cdot T)_{|_I}\rangle \\
&=&\langle m,(u_0u) \cdot T_{|_I}\rangle \\
&=&(u_0u)(x)\langle m,T_{|_I}\rangle \\
&=&u(x)(u_0(x)\langle m,T_{|_I}\rangle )\\
&=&u(x)\langle m,u_0\cdot (T_{|_I})\rangle \\
&=&u(x)\langle M,T\rangle \\
\end{eqnarray*}

It follows that $M\in TIM _{\cA(G)}(\cA(G)^*, \phi_x )$. Finally, if $T\in I^*$ and if $T_1\in \cA(G)$ with $T_{1_{|_I}}=T$, then 
\begin{eqnarray*}
\langle \hat{M},T\rangle &=&\langle M,T_1\rangle \\
&=&\langle m, T_{1_{|_I}}\rangle \\
&=& \langle m,T\rangle \\
\end{eqnarray*} 
Hence $\Gamma(M)=\hat{M}=m$ and $\Gamma$ is surjective.

\end{pf} 

The following result follows immediately from Corollary \ref{unique}. 

\begin{cor}\label{uniqueI} Let $\cA(G)$ be $A(G)$, $A_{cb}(G)$ or $A_M(G)$. Let $I$ be a closed ideal in $\cA(G)$. Assume that $x\not \in Z(I)$. If $G$ is a non-discrete group, then 
\[|   TIM _{I}(I^*, \phi_{{x}_{|_{I}} })|=2^{2^{b(G)}}.\]
In particular, $I^*$ admits a unique topological invariant mean if and 
only if $G$ is discrete.  

\end{cor}

\begin{lem} \label{restrict} Let $\cA$ be a commutative Banach algebra with maximal ideal space $\Delta(\cA)$. Let $X$ be a closed submodule of 
$\cA^*$ containing $\phi \in \Delta(\cA)$. Let $M \in TIM_{\cA}(\cA, \phi)$. Let $m=M_{|_X}$ be the restriction of $M$ to $X$. Then 
$m \in TIM_{\cA}(X, \phi)$. 

In particular, if we let $\cA(G)$ be one  of the algebras $A(G)$,  $A_{cb}(G)$ or $A_M(G)$,  $I$ be a closed ideal in $\cA(G)$ with 
 $x\not \in Z(I)$ and $X$ is any of $UCB(I)$, $WAP(I)$ or $AP(I)$, we have 
that $\phi_{x_{|_I}}\in X$ and hence that for any $M \in    TIM _{I}(I^*, \phi_{{x}_{|_{I}} })$, if $m=M_{|_X}$ we get that 
$m \in TIM _{I}(X, \phi_{{x}_{|_{I}} })$. 

\end{lem} 

\begin{pf} Let $M \in TIM_{\cA}(\cA, \phi)$. Let $m=M_{|_X}$ be the restriction of $M$ to $X$. We have that  
\[1=\langle M,\phi\rangle = \langle m,\phi\rangle =\|m\|_{X^*}\]
and that if $u\in \cA$ and $T\in X$, then 
\[\langle m, u\cdot T \rangle =\langle M, u\cdot T \rangle =\langle \phi, u\rangle \langle M, T\rangle  = \langle \phi, u\rangle \langle m, T\rangle .\]
Hence $m \in TIM_{\cA}(X, \phi)$.

Let $\cA(G)$ be one  of the algebras $A(G)$,  $A_{cb}(G)$ or $A_M(G)$ and   $I$  a closed ideal in $\cA(G)$ with 
 $x\not \in Z(I)$ . Let $X=UCB(I)$. We can choose an open neighborhood $U$ of $X$ with $U\cap Z(I)=\emptyset$ and a $u_0\in \cA(G)\cap C_{00}(G)$ such that $supp(u_0)\subseteq U$ and $u_0(x)=1$. Then $u_0\in I$. Moreover, if $v \in I$, 

\[\langle u_0\cdot \phi_{x_{|_I}}, v\rangle =\langle \phi_{x_{|_I}}, u_0v\rangle =u_0(x)v(x)=\langle \phi_{x_{|_I}}, u_0v\rangle .\]
Hence  $\phi_{x_{|_I}}=u_0\cdot \phi_{x_{|_I}} \in UCB(I)$. 

To see that $\phi_{x_{|_I}}\in AP(I)$ note that $\{u(x)|\|u\|_{I}\leq 1\} = \{\lambda \in \mathbb{C}| \mid \lambda \mid \leq 1\}$ is compact and hence 
\[ \{u\cdot \phi_{x_{|_I}}| \|u\|_{I}\leq 1\}= \{\lambda \phi_{x_{|_I}} \in \mathbb{C}| \mid \lambda \mid \leq 1\}\]
is compact in $I^*$ so $\phi_{x_{|_I}} \in AP(I)$. As $AP(I)\subseteq WAP(I)$ we also have that $\phi_{x_{|_I}} \in WAP(I)$.

\end{pf}

\begin{thm} Let $\cA(G)$ be one  of the algebras $A(G)$,  $A_{cb}(G)$ or $A_M(G)$. Let $I$ be a closed ideal in $\cA(G)$. Assume that $x\not \in Z(I)$. 
The restriction map $R:TIM_{I}(I^*, \phi_{x_{|_I}}))\rightarrow TIM_I(UCB(I), \phi_{x_{|_I}}))$ is a bijection.
In particular, if $G$ is non-discrete, then  
\[\mid TIM_{\cA(G)}(\cA(G)^*, \phi_x)\mid =\mid TIM_{\cA(G)}( UCB(\cA(G), \phi_x) \mid =2^{2^{b(G)}}.\]

\end{thm}

\medskip
\begin{pf}
Choose an open neighbourhood $U$ of $x$ with $U\cap Z(I)=\emptyset$ and a $u_0 \in \cA(G)\cap C_{00}(G)$ with $supp(u_0)\subseteq U$, and $u_0(x)=1$. Then $u_0 \in I$ and \[\phi_{x_{|_I}} =u_0 \cdot  \phi_{x_{|_I}} \in UCB(I).\]

Next we let $M \in TIM_{I}(I^*, \phi_{x_{|_i}}))$. Let $m=R(M)$. 
It follows from Lemma \ref{restrict} that $m\in TIM_I(UCB(I), \phi_{x_{|_I}}))$.

If
$M_{1},M_{2}\in TIM_{I}(I^*, \phi_{x_{|_I}})$ with $M_{1}\not=M_{2}$, 
then there exists a $T\in I^* $ for which 
\[\langle M_{1},T\rangle\not=\langle M_{2},T\rangle.\]

Choose an open neighbourhood $U$ of $x$ with $U\cap Z(I)=\emptyset$ and a $u_0 \in \cA(G)$ with $supp(u_0)\subseteq U$, and $u_0(x)=1$, then $u_0\cdot T \in UCB(I)$ with 
\[\langle M_{1},u_{0}\cdot T\rangle=\langle M_{1},T\rangle\not=\langle M_{2},T\rangle
=\langle M_{2},u_{0}\cdot T\rangle.\]
This shows that $R(M_{1})\not=R(M_{2})$ and hence $R$ is injective.

Next, let $m\in TIM_I(UCB(I), \phi_{x_{|_I}}))$. Pick a $u_{0}\in I$ with 
$\parallel u_{0}\parallel _{A(G)}=1=u_{0}(x)$. Define $M\in I^{**}$ by
\[\langle M,T\rangle =\langle m,u_{0}\cdot T\rangle, \, \; \; \; \; T\in I^*.\]
Since $u_{0}(x) = 1$, it follows that 
\[\langle M,\phi_{x_{|_I}}\rangle= \langle m,u_{0}\cdot \phi_{x_{|_I}}\rangle= \langle m,\phi_{x_{|_I}}\rangle=1.\]
From this and the fact that $\parallel u_{0}\parallel _{A(G)}=1$, we get that 
$\parallel M\parallel =1$. 

Next, if $v\in I, T\in I^*$ , then 
 \[ \langle M,v\cdot T\rangle = \langle m,u_{0}\cdot(v\cdot  T)\rangle 
 = \langle m,v\cdot(u_{0}\cdot  T)\rangle 
 = v(x) \langle m,u_{0}\cdot  T\rangle 
 = v(x) \langle M,  T\rangle. \]
 This shows that $M\in   TIM_{I}(I^*, \phi_{x_{|_I}})$. 

Finally, if $T\in UCB(I))$, then 
\[\langle M, T\rangle=\langle m, u_{0}\cdot T\rangle=\langle m, T\rangle\]
since $m\in TIM_I(UCB(I), \phi_{x_{|_I}})$. Therefore, $R(M)=m$
and $R$ is surjective. 
 
\end{pf}

\begin{rem} In the proof of the previous theorem we were able to explicitly show how each $m\in TIM_I(UCB(I), \phi_{x_{|_I}})$ extends to an element $M\in   TIM_{I}(I^*, \phi_{x_{|_I}})$. The next proposition  shows that such extensions hold in greater generality.

\end{rem}
We need the following lemma.

\begin{lem}\label{synth}  Let $\cA(G)$ be one  of the algebras $A(G)$,  $A_{cb}(G)$ or $A_M(G)$. Let $I$ be a closed ideal in $\cA(G)$ with $Z(I)$ being a set of spectral synthesis for $\cA(G)$. Assume that $x\not \in Z(I)$. Then $\{x\}$ is a set of spectral synthesis for $I$.

\begin{pf} Let $u\in I$ be such that $u(x)=0$. Let $E=Z(I)$. Let $\epsilon > 0$. Since $E$ is a set of spectral synthesis for $\cA(G)$, we can find $w\in \cA(G)\cap C_{00}(G)$ which is such  $K=supp(w)\cap E=\emptyset$ and 
\[\|u-w\|_{\cA(G)} <\frac{\epsilon}{2}.\]
Since $u(x)=0$, this means that $|w(x)|<\frac{\epsilon}{2}$.

Next we choose neighbourhoods $V_1$ and $V_2$ of $x$ with compact closure  disjoint from $E$ with  $V_1\subseteq V_2$. Then choose a $v\in \cA(G)$ such that $v(x)=0$ if $x\not \in V_2$, $v(y)=w(x)$ on $V_1$ and $\|v\|_{\cA(G)}=\mid w(x)\mid$. 
Then $w-v\in I=I(E)$ has compact support $K_1$ with 
$K_1\cap (E\cup \{x\})=\emptyset$ and 
\[\|u-(w-v)\|_{\cA(G)}\leq \|u-w\|_{\cA(G)} + \|v\|_{\cA(G)} < \frac{\epsilon}{2} + \frac{\epsilon}{2}= \epsilon.\]
\end{pf} 
\end{lem} 

\begin{prop} Let $\cA(G)$ be one  of the algebras $A(G)$,  $A_{cb}(G)$ or $A_M(G)$. Let $I$ be a closed ideal in $\cA(G)$ with $Z(I)$ being a set of spectral synthesis for $\cA(G)$. Assume that $x\not \in Z(I)$. Assume also that $Y\subset X$ are two closed submodules of $I^*$
each containing $\phi_{x_{|_I}}$. Let $m\in TIM_I(Y,\phi_{x_{|_I}})$. Then there exists some $M\in TIM_I(X,\phi_{x_{|_I}})$ such that $M_{|_Y}=m$. 

\end{prop}

\begin{pf} Let $m\in TIM_I(Y,\phi_{x_{|_I}})$. By the Hahn-Banach Theorem we can find a $\Psi \in X^*$ so that 
\[\|\Psi\|_{X^*}=\|m\|_{Y^*}=1.\]

Next we let 
\[S=\{u \in I\mid \|u\|_{\cA(G)}=u(x)=1\}.\]
Then $S$ is a commutative semigroup under pointwise multiplication and is hence amenable. Let $\Phi\in \ell_{\infty}(S)^*$ be an invariant mean.  Then for each $T\in X$ we define $f_T:S\to \mathbb{C}$ by 
\[f_T(u)=\langle \Psi, u \cdot T\rangle .\]
It follows that $f_T\in \ell_{\infty}(S)$ as $|f_T(u)|\leq \|T\|_{X}$ for each $u\in S$. Moreover, if $v\in S$, then  
\[f_{v\cdot T}(u)= \langle \Psi, u \cdot (v\cdot T)\rangle =\langle \Psi, vu \cdot T\rangle =L_v(f_T)(u) \]
where $L_v$ is the left translation operator on $\ell_{\infty}(S)$. 

Next, we let 
\[\langle M,T\rangle  = \langle \Phi, f_T\rangle \]
for each $T\in X$. Note that $f_{\alpha T_1+\beta T_2}=\alpha f_{T_1} + \beta f_{T_2}$ and $\mid\langle M,T\rangle \mid \leq \|T\|_{X}$ so that in fact $M\in X^*$ with $\|M\|_{X^*}\leq 1$. 

Now if $T=\phi_{x_{|_I}}$, then 
\[f_T(u)=\langle \Psi, u \cdot \phi_{x_{|_I}}\rangle  =1\] for every 
$u\in S$ and hence $ \langle M, \phi_{x_{|_I}}\rangle =1$ and $\|M\|_{X*}= 1$.

Finally, since $\Phi$ is a left invariant mean on $\ell_{\infty}(S)^*$ we have that if $u\in S$, then 
\[\langle M,u\cdot T\rangle =\langle \Phi, f_{u\cdot T}\rangle =\langle \Phi, L_u(f_T)\rangle =\langle M,T\rangle .\]

Finally, we must show that $\langle M,u\cdot T\rangle =u(x) \langle M,T\rangle $ for all $u\in I$, and $T\in X$.  

First we will show that if $v\in I$ and if there exists a neighborhood $U$ of $x$ such that $v(y)=1$ on $U$, then $\langle M,v\cdot T\rangle =\langle M,T\rangle $ for all $T\in X$. To see why this is the case, we choose $u\in I $ such that $u(x)=1=\|u\|_{\cA(G)}$ with $u(y)=0$ if $y\not \in U$. Then $uv=u$, hence 
\[\langle M, v\cdot T \rangle = \langle M, u\cdot(v\cdot T) \rangle =\langle M, uv\cdot T\rangle =\langle M, u\cdot T\rangle =\langle M,T\rangle .  \]

Next assume that $v\in I$ satisfies $v(y)=0$ on some neighborhood $U$ of $x$. Since $1\in B(G)$, the function $1-u$ is a multiplier of $I$, that is $(1-v)w\in I$ for every $w \in I$ and hence $(1-v)\cdot T \in X$.
We note that $1-v(y)=1$ on U. Again choosing $u\in I $ such that $u(x)=1=\|u\|_{\cA(G)}$ with $u(y)=0$ if $y\not \in U$. Once more we get that 
\[(\langle M,  T \rangle -\langle M, v\cdot T \rangle )=\langle M, (1-v)\cdot T \rangle = \langle M, u\cdot((1-v)\cdot T) \rangle  = \langle M, u\cdot T\rangle =\langle M,T\rangle .  \]
Hence $\langle M, v\cdot T \rangle =0$.

Now let $u\in I$ be such that $u(x)=0$. By Lemma \ref{synth}, $\{x\}$ is a set of spectral synthesis for $I$. It follows that we can find a sequence of functions $\{w_n\}\subset I$ and as sequence $\{U_n\}$ of open neighbourhoods of $x$ such that each $w_n$ has compact support, $w_n(y)=0$ for all $y\in U_n $ and $\lim\limits_{n\to \infty} \|u-w_n\|_{\cA(G)}=0$. In particular, from what we have seen above, 
$\langle M, w_n\cdot T\rangle =0$ and hence 
\[\langle M, u\cdot T\rangle = \langle M, (u-w_n)\cdot T\rangle \]
Moreover, 
\[\lim\limits_{n\to \infty} \|u\cdot T - w_n\cdot T\|_{I^*}=0,\]
as hence we have that $\langle M,u \rangle=0$. 
Finally, choose $u\in I$ be such that $u(x)=1$. Once more choose 
$v\in I$ so that $v=1$ on a neighbourhood $U$ of $x$. Then $u-v(x)=0$ and 
\[\langle M,(u-v)\cdot T\rangle =0 \]
which means that 
\[\langle M,u\cdot T\rangle =\langle M,v\cdot T\rangle =v(x)\langle M,T\rangle =u(x)\langle M,T\rangle .\]

For here, if $u(x)\not =0$, let $w=\frac{1}{u(x)} u $. Then 
\[\langle M,u\cdot T\rangle = \langle M,u(x)( \frac{1}{u(x)} u\cdot T)\rangle =u(x) \langle M,( \frac{1}{u(x)} u\cdot T)\rangle =u(x)\langle M,T\rangle .\]
This shows that $M\in TIM_I(X,\phi_{x_{|_I}})$. 

Finally, if $T\in Y$, then 
\[f_T(u)=\langle \Psi, u \cdot T\rangle = \langle m, u \cdot T\rangle =\langle m,T\rangle \]
for all $u\in S$. In particular, 
$\langle M,T\rangle = \langle \Phi,f_T\rangle =\langle m,T\rangle $ so 
$M_{|_Y}=m$ as desired.

\end{pf}

\begin{thm} Let $\cA(G)$ be one  of the algebras $A(G)$,  $A_{cb}(G)$ or $A_M(G)$. Let $I$ be a closed ideal in $\cA(G)$. Assume that $x\not \in Z(I)$. Then $WAP(I)$ has a unique topologically invariant mean at $\phi_{x_{|_I}}$.

\end{thm} 

\begin{pf} The fact that $TIM_{I}(WAP(I), \phi_{x_{|_i}}))\not =\emptyset$ follows immediately from the observation that every TIM on $I^*$ restricts to a TIM on $WAP(I)$ and we know that $TIM_{I}(I^*, \phi_{x_{|_i}}))\not =\emptyset$.

Next, we let $n,m\in TIM_{I}(WAP(I), \phi_{x_{|_i}}))$.  Then there exits a net $\{u_{\alpha}\}_{\alpha\in \Omega}\subseteq I$ so that
\[m=\lim\limits_{\alpha\in \Omega}\{u_{\alpha}\}\]
in the weak$^*$ topology on $I^{**}$. Hence, for each $T\in I^*$, 
\[\langle n \odot m,T \rangle = \lim\limits_{\alpha}\langle n\odot u_{\alpha}, T\rangle =\lim\limits_{\alpha}\langle n, u_{\alpha}\cdot T\rangle =\lim\limits_{\alpha}u_{\alpha}(x)\langle n, T\rangle .\]
But we also know that 
\[1=\langle n, \phi_{x_{|_i}}\rangle =\lim\limits_{\alpha}\langle  u_{\alpha},\phi_{x_{|_i}}\rangle =\lim\limits_{\alpha} u_{\alpha}(x).\]
Hence $n\odot m= n$. But we also know that $n\odot m=m\odot n=m$. Hence, $n=m$ and the TIM is unique. 
\end{pf} 

\begin{cor} Let $\cA(G)$ be one  of the algebras $A(G)$,  $A_{cb}(G)$ or $A_M(G)$. Let $I$ be a closed ideal in $\cA(G)$. If $UCB(I)\subseteq WAP(I)$, then $G$ is discrete.

\end{cor} 

\begin{pf} Assume that $G$ is non discrete and that $x\not \in Z(I)$. Let $M_1, M_2 \in TIM_{I}(I^*, \phi_{x_{|_i}}))$ with $M_1\not = M_2$ , Then 
$m_1=M_{1_{|_{WAP(I)}}}, m_2=M_{2_{|_{WAP(I)}}} \in TIM_{I}(WAP(I), \phi_{x_{|_i}}))$. Let $T\in I^*$ be such that 
$\langle M_1, T\rangle \not =\langle M_2,T\rangle $ and choose $u\in I$ so that $u(x)=1$. Then $u\cdot T \in UCB(I) \subseteq WAP(I)$, and 
\[\langle m_1, u\cdot T\rangle = \langle M_1, u\cdot T\rangle =\langle M_1,T\rangle \not = \langle M_2,T\rangle =\langle M_2, u\cdot T\rangle =\langle m_2, u\cdot T\rangle \]
which contradicts the uniqueness of the TIM on $WAP(I)$.

\end{pf}

\section{Arens Regularity of Ideals in \texorpdfstring{$A(G)$, $A_{cb}(G)$ and $A_M(G)$}{} } 

In this section, we will apply what we know about topologically invariant means to questions concerning the possible Arens regularity of ideals in $A(G)$, $A_{cb}(G)$, and $A_M(G)$. The key observation is the following which improves on \cite[Corollary 3.13]{For2}: 

\begin{thm} \label{discrete} Let $\cA(G)$ be any of the algebras $A(G)$, $A_{cb}(G)$ or $A_M(G)$. Let I be a non-zero closed ideal in $\cA(G)$. If $I$ is Arens regular, then $G$ is discrete. 

\end{thm} 

\begin{pf} If $I$ is Arens regular, then $I^*=WAP(I)$. Hence $I^*$ has a unique topologically invariant mean. However, by Corollary \ref{uniqueI}, this implies that $G$ must be discrete. 

\end{pf} 

The following corollary is immediate. See also \cite[Theorem 3.2]{For2} and \cite[Corollary 3.9]{ForMiao}. 
 \begin{cor} Let $\cA(G)$ be any of the algebras $A(G)$, $A_{cb}(G)$ or $A_M(G)$. If $\cA(G)$ is Arens regular, then $G$ is discrete.  

\end{cor}

\begin{cor} Let $G$ be non-discrete. If $\cA(G)$ is one of $A(G)$, $A_{cb}(G)$ or$ A_M(G)$, then $\cA(G)$ has no non-zero reflexive closed ideal. 
    
\end{cor}

\begin{rem} Granirer \cite[Theorem 5]{Gran} has shown that every infinite discrete group contains an infinite  set $E\subset G$ such that the ideal $I_{A(G)}(G\setminus E)$ is isomorphic to $\ell_2$. In particular, this ideal is reflexive and hence Arens regular.  However, if we ask that $I$ also has a bounded approximate identity, then at least for the Fourier algebra such an ideal can only be Arens regular if it is finite-dimensional. 
    
\end{rem}

\begin{lem}
Let $\cA(G)$ be any of the algebras $A(G)$, $A_{cb}(G)$ or $A_M(G)$. Let $H$ be a subgroup of $G$. If $\cA(G)$ is Arens regular, then so is $\cA(H)$. In particular, if $H$ is amenable, then $H$ is finite.   
\end{lem} 

\begin{pf} As $G$ is discrete, $H$ is open in $G$. In this case, the restriction map $R:\cA(G)\to \cA(H)$ is a contractive homomorphism that is also surjective. As such $\cA(H)$ is Arens regular. 

The last statement is simply \cite[Proposition 3.3]{Lau2}.

\end{pf}

\begin{defn} Let $\mathcal{R}(G)$, the \it{coset ring of} $G$,   denote the Boolean ring of sets generated by cosets of subgroups of $G$.  A subset $E$ of $G$ is in  $\mathcal{R}(G)$ if and only if \[E=\bigcup\limits_{i=1}^n ( x_i H_i \setminus \bigcup\limits_{j=1}^{m_i} b_{i,j}K_{i,j}), \]
where $H_i$ is a  subgroup of $G$, $x_i\in G$, $K_{i,j}$ is a subgroup of $H_i$, and $b_{i,j}\in K_{i,j}$. 

By $\mathcal{R}_a(G)$, the \it{amenable coset ring of} $G$, we will mean all sets of the form \[E=\bigcup\limits_{i=1}^n ( x_i H_i \setminus \bigcup\limits_{j=1}^{m_i} b_{i,j}K_{i,j}), \]
where $H_i$ is an amenable subgroup of $G$, $x_i\in G$, $K_{i,j}$ is a subgroup of $H_i$, and $b_{i,j}\in K_{i,j}$. 
  
\end{defn}

\begin{thm}\label{Amen1} Let $I$ be a closed ideal of $A(G)$ with a bounded approximate identity that is Arens regular. Then $I$ is finite-dimensional. 
\end{thm} 

\begin{pf} Assume that $I\subseteq A(G)$ is Arens regular. By Theorem \ref{discrete} $G$ must be discrete. If $I$ has a bounded approximate identity, then since $A(G)$ is weakly sequentially complete (ref), $I$ must be unital. It follows that $1_{G\setminus Z(I)}\in I$. In particular, $G\setminus Z(I)$ must be compact and hence finite. This shows that $I$ is finite-dimensional.  

\end{pf} 

\begin{rem} The fact that $A(G)$ is  weakly sequentially complete was crucial in establishing the previous theorem. Unfortunately, we do not know whether or not either or both of $A_{cb}(G)$ or $A_M(G)$ would be weakly sequentially complete. 

\end{rem} 

For the remainder of this section, we will assume that $G$ is a discrete group. 

\begin{lem}\label{Amensubgroup} Let $\cA(G)$ be any of the algebras $A(G)$, $A_{cb}(G)$ or $A_M(G)$. Let $H$ be a proper amenable subgroup of $G$. If $I_{\cA(G)}(H)$ is Arens regular, then $H$ is finite and $\cA(G)$ is also Arens regular. 

\end{lem} 
\begin{pf} Since $H$ is proper, there exists an $x\in G\setminus H$. The ideal $I_{\cA(G)}(xH)$ is isometrically isomorphic to $I_{\cA(G)}(H)$ and hence is also Arens regular. 

If $u \in \cA(H)$, then the function $u^{\circ}$ defined by $u^{\circ}(y)=u(y)$ if $y\in H$ and $u^{\circ}(y)=0$ if $y\in G\setminus H$ is in $\cA(G)$. Now if $R:\cA(G)\to \cA(H)$ is the restriction map, then $R$ is contractive homomorphism that maps $I_{\cA(G)}(xH)$ onto $\cA(H)$. In particular, $\cA(H)=A(G)$ is also Arens regular. It follows that $H$ is finite. 

Let $\mathcal{B}$ be the algebra $1_H\cA(G)\oplus_1 1_{G\setminus H} \cA(G)$. Then $B$ is a commutative Banach algebra and the mapping $\Gamma:\cA(G)\to \mathcal{B}$ given by $\Gamma(u) = (1_Hu, 1_{G\setminus H}u)$ is a continuous isomorphism  that maps
$I_{\cA(G)}(H)$ isometrically onto the ideal $(I_{\cA(G)}(H), 0)$ in $\mathcal{B}$. Since $1_{G\setminus H})\cA(G)$ is finite-dimensional, it is Arens regular. We get that 
\[(1_H\cA(G)\oplus_1 1_{G\setminus H} \cA(G))^{**}=(1_H\cA(G))^{**}\oplus_1 (1_{G\setminus H} \cA(G))^{**}\]
which is commutative since each of its components is commutative. Hence 
$1_H\cA(G)\oplus_1 1_{G\setminus H} \cA(G)$ is Arens regular, and so is $\cA(G)$
\end{pf} 

\begin{thm} Let $\cA(G)$ be any of the algebras $A(G)$, $A_{cb}(G)$ or $A_M(G)$. Let \[E=\bigcup\limits_{i=1}^n ( x_i H_i \setminus \bigcup\limits_{j=1}^{m_i} b_{i,j}K_{i,j}), \]
where $H_i$ is an amenable subgroup of $G$, $x_i\in G$, $K_{i,j}$ is a subgroup of $H_i$, and $b_{i,j}\in K_{i,j}$. If $I_{\cA(G)}(E)$ is non-zero and Arens regular, then either $E$ is finite and $\cA(G)$ is also Arens regular, or $G$ is amenable and $I(E)$ is finite-dimensional. 

\end{thm} 
\begin{pf}

We begin by first assuming that $E=\bigcup\limits_{i=1}^n  x_i H_i $. In this case, we will prove the conclusion by induction on $n$. That is we let $P(n)$ be the statement that if $E=\bigcup\limits_{i=1}^n  x_i H_i $ is a proper subset of $G$ and if $I_{\cA(G)}(E)$ is Arens regular, then $E$ is finite and $\cA(G)$ is Arens regular. 

If $n=1$, $E=xH$ where $H$ is a proper amenable subgroup. Since $I(H)$ is isometrically isomorphic to $I(xH)$, Lemma \ref{Amensubgroup} shows that $E$ is finite and $\cA(G)$ is Arens regular.

Assume that $P(n)$ is true for all $n\leq k$. Let $E=\bigcup\limits_{i=1}^{k+1}  x_i H_i $ where each $H_i$ is an amenable subgroup of $G$. By translating if necessary we can assume that $x_{k+1}=e$. If $H_{k+1}\subseteq \bigcup\limits_{i=1}^{k}  x_i H_i $, then we have $E=\bigcup\limits_{i=1}^{k}  x_i H_i $ and we are done. So we may assume that 
\[F=H_{k+1}\setminus (\bigcup\limits_{i=1}^{k}  x_i H_i)\not = \emptyset.\]

Note that $H_{k+1}\setminus F \in \mathcal{R}(H_{k+1}) $  and 
\[I_{\cA(H)}(H_{k+1}\setminus F)=I_{\cA(G)}(E)_{|_{H_{k+1}}}.\]
In particular, since the restriction map is a homomorphism, $I_{\cA(H)}(H_{k+1}\setminus F)$ is Arens regular. But as $H_{k+1}\setminus F \in \mathcal{R}(H_{k+1}) $ and $H_{k+1}$ is amenable, we have that $\cA(H)=A(H)$ and $I_{\cA(H)}(H_{k+1}\setminus F)$ has a bounded approximate identity. It then follows from Theorem \ref{Amen1} that $F$ is finite. 

Next we observe that $E$ is the disjoint union of $\bigcup\limits_{i=1}^{k}  x_i H_i $ and the finite set $F$. But as $F$ is finite we can proceed in a manner similar to that of the proof of Lemma \ref{Amensubgroup} to conclude that $I_{\cA(G)}(\bigcup\limits_{i=1}^{k}  x_i H_i)$ is also Arens regular. From here the induction hypothesis tells us that $\bigcup\limits_{i=1}^{k}  x_i H_i$ is finite. And as $F=H_{k+1}\setminus (\bigcup\limits_{i=1}^{k}  x_i H_i)$ is also finite, $H_{k+1}$ is finite. Hence $E$ is finite as well. 

If we assume that \[E=\bigcup\limits_{i=1}^n ( x_i H_i \setminus \bigcup\limits_{j=1}^{m_i} b_{i,j}K_{i,j}), \]
where $H_i$ is an amenable subgroup of $G$, $x_i\in G$, $K_{i,j}$ is a subgroup of $H_i$. We have two cases. The first is that $\bigcup\limits_{i=1}^n  x_i H_i \not = G$. If this is the case, then if $E_1=\bigcup\limits_{i=1}^n  x_i H_i $ then $E\subseteq E_1$ and hence the non-zero closed ideal $I_{\cA(G)}( E_1) $ is contained in the Arens regular ideal $I_{\cA(G)}( E)$  and is therefore also Arens regular. But we have seen above that this means that $E_1$ is finite. It follows that $E$ is also finite. As before, this would imply that $\cA(G)$ would also be Arens regular. 

Finally,  if we assume that $\bigcup\limits_{i=1}^n  x_i H_i = G$. Then by \cite[Corollary 3.3]{For1} one of the $H_i$'s has finite index in $G$. Since each $H_i$ is amenable, so is $G$. This means that we can express $G\backslash E$ as a disjoint union 
$\bigcup\limits_{l=1}^{m} F_l$ where each $F_l$ is a translate of an element of the  coset ring of one of the open amenable subgroups $K_{i,j}$. Moreover, this means that $I_{\cA(G)}(E)=I_{A(G)}(E)$ has a bounded approximate identity \cite[Theorem 3.20]{For1}. It now follows from Theorem \ref{Amen1} that this ideal is finite-dimensional.

\end{pf}

\vskip 0.4 true cm

\begin{center}{\textbf{Acknowledgments}}
\end{center}
The authors wish to congratulate Professor Anthony To-Ming Lau on his many accomplishments over his distinguished career. In addition, the first author wishes to express his deep gratitude to Professor Lau for his exceptional mentorship, his unwavering support, and for the exemplary role model he has provided for all of us who have had the great pleasure to work with him.  \\ \\
\vskip 0.4 true cm

%------------------------------------------------------------------------------------%

%
\vfill

{\footnotesize \par\noindent{\bf Brian Forrest}}\; \\ {Department of Pure Mathematics}, {University
of Waterloo} {Waterloo, Ontario, Canada, N2L 3G1}\\
{\tt Email: beforres@uwaterloo.ca}\\

{\footnotesize \par\noindent{\bf John Sawatzky}\; \\ {Department of Pure Mathematics}, {University
of Waterloo,} {Waterloo, Ontario, Canada, N2L 3G1}\\
{\tt Email: jmsawatzky@uwaterloo.ca}}\\

{\footnotesize \par\noindent{\bf Aasaimani Thamizhazhagan}\; \\ {Department of Pure Mathematics}, {University
of Waterloo, } {Waterloo, Ontario, Canada, N2L 3G1}\\
{\tt Email: athamizhazhagan@uwaterloo.ca}}

\end{document}